\documentclass[11pt]{article}
\usepackage{amsmath}
\usepackage{amsthm}
\usepackage{amssymb}
\usepackage{graphics}
\usepackage{bm}
\usepackage{setspace}
\usepackage[sort&compress]{natbib}
\usepackage{authblk}
\usepackage{float}
\usepackage{comment}
\usepackage{tikz}
\usepackage{tikz-cd}
\usepackage{amsfonts}
\usepackage[cal=bickham,calscaled=1.05]{mathalfa}
\usepackage{eucal}
\usepackage{latexsym}
\usepackage{mathrsfs}
\usepackage{mathtools}
\usepackage{hyphenat}
\usepackage{amscd}
\usepackage{enumerate}
\usepackage[pagebackref,hypertexnames=false, colorlinks, citecolor=blue, linkcolor=red, pdfstartview=FitB]{hyperref}
\usepackage[utf8]{inputenc}
\usetikzlibrary{positioning}
\usepackage{fullpage}
\usepackage{color}
\numberwithin{equation}{section}
\numberwithin{subsection}{section}
\newtheorem{thm}{Theorem}[section]

\newtheorem{prop}[thm]{Proposition}

\newtheorem*{thm*}{Theorem}
\theoremstyle{definition}

\newtheorem{remark}[thm]{Remark}
\newtheorem{eg}[thm]{Example}

\title{$L^p$-Fourier--Stieltjes Algebras of Certain Locally Compact Groups}
\author[]{Hamidreza Shirnejad\\[1em]
\normalsize University of Manitoba
}

\affil{\footnotesize}

\date{}

\newcommand{\keywords}[1]{%
    \par\medskip
    \noindent\textbf{Keywords:}\,#1
    \par\medskip
}

\begin{document}

\bibliographystyle{plain}
\maketitle
%%%%%%%%%%%%%%%%%%%%%%%%%%%%%%%%%%%%%%%%%%%%%%%%%%%%%%%%%%%%%%%%%%%%%%%%%%%%%%%%%%%%%%%%%%%%%%%%%%%%%%%%%%%%
\begin{abstract}
We study the relationship between the Brown--Guentner (BG) and Kaliszewski--Landstad--Quigg (KLQ) constructions of $L^p$-Fourier--Stieltjes algebras. For certain nonamenable locally compact groups satisfying the Kunze--Stein property with an exponential integrability condition and an ideal classification assumption, as well as a class of semidirect product groups $G = N \rtimes H$, under an appropriate stabilizer assumption, we show that $B_{L^p,\text{BG}}(G) = B_{L^p,\text{KLQ}}(G)$, for every $p \in [2, \infty)$. Motivated by the known classifications of the Fourier--Stieltjes ideals and their connection with the unitary dual, we also describe the closed ideals of $\widehat{G}$ for these semidirect product groups in terms of those of $\widehat{H}$.
\end{abstract}



%%%%%%%%%%%%%%%%%%%%%%%%%%%%%%%%%%%%%%%%%%%%%%%%%%%%%%%%%%%%%%%%%%
\section{Introduction}

Let $G$ be a locally compact group. Associated with $G$ are two canonical $C^*$-algebras: the full group $C^*$-algebra $C^*(G)$ and the reduced group $C^*$-algebra $C^*_{r}(G)$. The full group $C^*$-algebra is obtained by completing the
convolution algebra $C_c(G)$ with respect to the universal $C^*$-norm coming from all cyclic unitary representations
of $G$, while the reduced group $C^*$-algebra is defined using the left regular representation of $G$. These two $C^*$-algebras coincide exactly when $G$ is amenable
(see \cite{Br-Gn}). For nonamenable groups they are
distinct, and it is natural to investigate
$C^*$-completions of $C_c(G)$ which lie between
the full and reduced group $C^*$-algebras.

The dual space of the full group $C^*$-algebra $C^*(G)$ has a description in terms of the representation theory of $G$. Indeed, the dual space of $C^*(G)$ can be identified with the Fourier--Stieltjes algebra $B(G)$, which consists of all coefficient
functions of continuous unitary representations of $G$. Through this duality, weak$^*$-closed
$G$-invariant ideals of $B(G)$ correspond to quotients of $C^*(G)$ (see \cite{Kan-Lau}). 

$C^*$-algebras that lie strictly between $C^*(G)$ and $C^*_r(G)$ are called ``exotic`` group $C^*$-algebras. An important source of such $C^*$-algebras comes from integrability conditions on coefficient functions of continuous unitary representations of $G$. Given $p \in [1, \infty]$, a unitary representation of $G$ is called an $L^p$-representation if it admits a dense subspace whose coefficient functions belong
to $L^p(G)$. These $L^p$-representations
give rise to $C^*$-norms on $C_c(G)$ and hence to $C^*$-algebras that may lie strictly between the full and reduced group $C^*$-algebras when $G$ is nonamenable. From the perspective of the Fourier--Stieltjes algebra, these constructions produce certain weak$^*$-closed $G$-invariant ideals of $B(G)$, referred to as $L^p$-Fourier--Stieltjes algebras. Two approaches to defining $L^p$-Fourier--Stieltjes algebras have appeared in the literature. The first approach is based on the coefficient functions
of $L^p$-representations introduced by Brown and Guentner through the associated group $C^*$-completions for discrete groups
(see \cite{Br-Gn}). Subsequently, Wiersma studied the corresponding $L^p$-Fourier and $L^p$-Fourier--Stieltjes algebras for locally compact groups (see \cite{WRS-Fnc}). A second approach, considered by Kaliszewski, Landstad, and Quigg, takes the weak$^*$-closure of $B(G)\cap L^p(G)$ (see \cite{KLQ-Ex-NY}, \cite{KLQ-Oi-NY}). 

An important question is whether these two constructions coincide. The question is equivalent to determining whether the weak$^*$-closure of $\text{span}(P(G)\cap L^p(G))$, where $P(G)$ denotes the set of continuous positive definite functions on $G$, coincides with the weak$^*$-closure of $B(G)\cap L^p(G)$. Kaliszewski, Landstad,
and Quigg studied the notions of ordered and weakly ordered invariant ideals of $B(G)$ and their relationship with this question (see \cite{KLQ-Oi-NY}).  Forrest, Tanko, and Wiersma obtained partial results on the coincidence of these constructions. In particular, they provided a positive answer to the corresponding question for $G = \text{SL}(2,\mathbb{R})$ (see \cite{FTW}). However, the general question remains open.

For certain groups, $L^p$-Fourier--Stieltjes algebras have been studied in detail and shown to characterize the weak$^*$-closed $G$-invariant ideals of $B(G)$. Such classifications
have been obtained for the groups $\text{SL}(2,\mathcal{K})$, where
$\mathcal{K} = \mathbb{R}$, $\mathbb{C}$, or $\mathbb{Q}_p$, and for suitable noncompact closed subgroups of automorphism groups of homogeneous trees satisfying appropriate transitivity and independence conditions (see \cite{WRS-Fnc}, \cite{Sie-FS-arx}, \cite{HLS-Tr-IMRN}).

The purpose of this work is to further study the relationship between the
BG and KLQ $L^p$-Fourier--Stieltjes algebras for certain classes of nonamenable locally compact groups. Our approach is based on harmonic analytic and representation theoretic methods, with an emphasis on the structure of the unitary dual.

The assumptions considered in this work are motivated by well-studied examples, including the groups $\text{SL}(2, \mathcal{K})$, where $\mathcal{K} = \mathbb {R}$, $\mathbb {C}$, or $\mathbb {Q}_p$, and certain groups acting on homogeneous trees. Indeed, these groups satisfy the Kunze--Stein property and admit an exponential integrability condition associated with conditionally negative definite continuous functions. They also have the property that BG-$L^p$-Fourier--Stieltjes algebras characterize the weak$^*$-closed $G$-invariant ideals of $B(G)$. Together, these three properties provide a framework for establishing the coincidence result for nonamenable locally compact groups satisfying them. More precisely, we prove that
\begin{align*}
B_{L^p,\text{BG}}(G) = B_{L^p,\text{KLQ}}(G), \hspace{3mm} 2\leq p <\infty.
\end{align*}
This provides further positive examples of the coincidence problem discussed above.

Motivated by the semidirect product group $\mathbb {R}^2 \rtimes \text{SL}(2, \mathbb {R})$, we study the coincidence of the BG and KLQ $L^p$-Fourier--Stieltjes algebras for a class of semidirect product groups $G = N \rtimes H$, where $N$ is a noncompact abelian locally
compact group and $H$ is a nonamenable type I group. We assume that $G$ is second
countable, and the stabilizers of non-fixed characters
in $\widehat{N}$ are amenable and type I. The amenability assumptions on the stabilizers are restrictive. Indeed, for a local field $\mathcal{K}$, the stabilizers of nontrivial characters in the semidirect product group $\mathcal{K}^3 \rtimes \text{SL}(3, \mathcal{K})$ are nonamenable. Thus, this group falls outside the scope of our framework.

Under the mentioned assumptions, Mackey theory allows us to describe $\widehat{G}$ in terms of the $H$-orbits in $\widehat{N}$. Representations associated with non-fixed characters are induced from amenable stabilizers and are weakly contained in the left regular representation $\lambda_G$. In contrast, representations associated with fixed characters have the form $\chi \times \pi$, where $\chi \in \widehat{N}$ is a character fixed by $H$, and $\pi \in \widehat{H}$, and their nonzero coefficient functions do not vanish at infinity. To show the coincidence of the two $L^p$-Fourier--Stieltjes constructions for these groups, we consider the Rajchman algebra $B_0(G) = B(G) \cap C_0(G)$. Using the structure of the unitary dual and type I decomposition theory, we prove that 
\begin{align*}
\overline{B_0(G)}^{\sigma(B(G),C^*(G))} = B_r(G),
\end{align*}
provided either the set of $H$-fixed characters in $\widehat{N}$ is countable or $B_0(G)|_H\subset B_r(H)$. This weak$^*$-closure identity is the key to obtain
\begin{align*}
B_{L^p,\text{BG}}(G) = B_r(G) = B_{L^p,\text{KLQ}}(G),
\hspace{3mm} 2\leq p<\infty.
\end{align*}
Thus, for this class of semidirect product groups, the two BG and KLQ constructions coincide and do not produce intermediate ideals between $B_r(G)$ and $B(G)$.

The classification of Fourier--Stieltjes ideals for $\text{SL}(2, \mathcal{K})$, where $\mathcal{K} = \mathbb{R}, \mathbb{C},$ or $\mathbb{Q}_p$, and certain groups acting on homogeneous trees provides further motivation for our study of semidirect product groups. For those groups, the $L^p$-Fourier--Stieltjes algebras describe the weak$^*$-closed $G$-invariant ideals of $B(G)$, which correspond to the closed ideals of the unitary dual $\widehat{G}$. This leads to the question of whether a similar classification can be found for the semidirect product groups considered here. We address this question by describing the closed ideals of $\widehat{G}$ in terms of those of $\widehat{H}$. More precisely, every such ideal is obtained by adjoining the reduced dual $\widehat{G}_r$ to representations of the form $\chi \times \pi$, where $\chi \in \widehat{N}$ is a $H$-fixed character and $\pi$ belongs to a closed ideal of $\widehat{H}$.

%%%%%%%%%%%%%%%%%%%%%%%%%%%%%%%%%%%%%%%%%%%%%%%%%%%%%%%%%%%%%%%%%%%%
\section{Preliminaries}

In this section, we recall some standard definitions and results needed for our work.

\subsection*{$L^p$-Fourier and $L^p$-Fourier--Stieltjes Algebras}

Let $G$ be a locally compact group and $p \in [1, \infty]$ be given. We say a unitary representation $\pi: G \to \mathcal{U}(\mathcal{H}_\pi)$ is a \emph{$\text{BG}$-$L^p$-representation} of $G$ when $\mathcal{H}_\pi$ admits a dense subspace $\mathcal{H}_0$ such that for each $\xi, \eta \in \mathcal{H}_0$, the coefficient function $\pi_{\xi, \eta}: G \to \mathbb{C}$, $g \mapsto \langle\pi(g)\xi, \eta \rangle$ lies in $L^p(G)$ (see \cite{Br-Gn}). Since coefficient functions are bounded, $\pi$ is a $\text{BG}\text{-}L^q$-representation for each $q > p$. We note that the left regular representation $\lambda_G: G \to \mathcal{U}(L^2(G))$ is a $\text{BG}\text{-}L^p$-representation of $G$ for each $ p \in [1, \infty]$ since $\lambda_{\xi, \eta} \in C_c(G)$ for every $\xi, \eta \in C_c(G)$. The \emph{BG-$L^p$-Fourier} and \emph{$\text{BG}$-$L^p$-Fourier--Stieltjes algebras} of $G$ are defined as 
\begin{align*}
    A_{L^p, \text{BG}}(G) := \overline{\text{span}\{\pi_{\xi , \eta}: \pi \hspace{2mm}\text{is an} \hspace{2mm} \text{BG}\text{-}L^p\text{-representation of }G, \,\xi, \eta \in \mathcal{H}_\pi\}}^{||\cdot||_{B(G)}} \subset B(G),
\end{align*}
and
\begin{align*}
    B_{L^p, \text{BG}}(G) := \overline{A_{L^p, \text{BG}}(G)}^{\sigma(B(G), C^*(G))} \subset B(G).
\end{align*}
It is proved that for each $p\in [1, \infty]$, $A_{L^p,\text{BG}}(G) = \overline{\text{span}(P(G) \cap L^p(G))}^{||\cdot||_{B(G)}}$, where $P(G)$ denotes the space of continuous positive definite functions on $G$ (see \cite{WRS-Fnc}).

For a representation $\pi:G \to \mathcal{U}(\mathcal{H}_\pi)$ we let
\begin{align*}
    A_\pi(G) := \overline{\text{span}\{\pi_{\xi, \eta}: \xi, \eta \in \mathcal{H}_\pi\}}^{||\cdot||_{B(G)}}, \quad \text{and} \quad B_\pi(G) := \overline{A_\pi(G)}^{\sigma(B(G), C^*(G))}.
\end{align*}

Kaliszewski, Landstad, and Quigg proposed another construction of $L^p$-Fourier and $L^p$-Fourier--Stieltjes algebras to be defined by 
\begin{align*}
    A_{L^p, \text{KLQ}}(G) := \overline{B(G) \cap L^{p}(G)}^{||\cdot||_{B(G)}}, \hspace{3mm} \text{and} \hspace{3mm} B_{L^p, \text{KLQ}}(G) := \overline{A_{L^p, \text{KLQ}}(G)}^{\sigma(B(G), C^*(G))},
\end{align*}
 and they conjectured that these two constructions coincide with the BG-constructions; however, there is very little evidence that this is true in general (see \cite{KLQ-Ex-NY}, \cite{KLQ-Oi-NY}). We say a unitary representation $\pi: G \to \mathcal{U}(\mathcal{H}_\pi)$ is a \emph{$\text{KLQ}\text{-}L^p$-representation} of $G$ if there exist dense subspaces $\mathcal{H}_1, \mathcal{H}_2 \subset \mathcal{H}_\pi$ such that $\pi_{\xi, \eta} \in L^p(G)$ for all $\xi \in \mathcal{H}_1$ and $\eta \in \mathcal{H}_2$. Also, we note that every $\text{BG}\text{-}L^p$-representation of $G$ is a $\text{KLQ}\text{-}L^p$-representation which implies that $B_{L^p, \text{BG}}(G) \subset B_{L^p, \text{KLQ}}(G)$ for all $p \in [1, \infty]$. It is shown that 
\begin{align*}
    A_{L^p, \text{KLQ}}(G) = \overline{\text{span} \{\pi_{\xi , \eta}: \pi \hspace{2mm}\text{is a}\hspace{2mm}\text{KLQ}\text{-}L^p\text{-representation of }G, \, \xi, \eta \in \mathcal{H}_\pi\}}^{||\cdot||_{B(G)}} \, \text{(see \cite{FTW})}.
\end{align*} 

We note that both $A_{L^p, \text{BG}}(G)$ and $A_{L^p, \text{KLQ}}(G)$ are contained in the \emph{Rajchman algebra} $B_0(G)$, where $B_0(G) := B(G) \cap C_0(G)$, for every $p \in [1, \infty)$. This follows from the left uniform continuity of coefficient functions and the left invariance of Haar measure.

\subsection*{Mackey Theory for Semidirect Product Groups}

Let $G = N \rtimes H$, where $N$ is an abelian locally compact group and $H$ is a locally compact group acting continuously on $N$. We write the group law on $G$ as
\begin{align*}
    (n_1, h_1)(n_2, h_2) = (n_1(h_1\cdot n_2),h_1 h_2), \hspace{3mm} n_1, n_2 \in N, \hspace{1mm} h_1, h_2 \in H.
\end{align*}
Identifying $N$ and $H$ with their canonical closed subgroups of $G$, we note that $N$ is a closed normal subgroup of $G$. The action of $H$ on $N$ induces a continuous action on the dual group $\widehat{N}$, given by
\begin{align*}
    (h \cdot \chi)(n) = \chi(h^{-1} \cdot n), \hspace{3mm} \chi \in \widehat{N}, \, h \in H, \,n \in N.
\end{align*}

For each $\chi \in \widehat{N}$, let 
\begin{align*}
  H_\chi := \{h \in H: h \cdot \chi = \chi\} \subset H \hspace{3mm} \text{and}  \hspace{3mm} H(\chi) :=\{h\cdot \chi: h \in H\} \subset \widehat{N}
\end{align*}
denote its \emph{stabilizer} and its \emph{orbit} under $H$, respectively. We also denote by $\widehat{N}/H$ the orbit space of $\widehat{N}$ under the action of $H$. A \emph{cross-section} of the $H$-orbits in $\widehat{N}$ is a subset $X \subset \widehat{N}$ meeting each orbit in exactly one point. Since $N$ acts trivially on $\widehat{N}$, the action of $G$ on $\widehat{N}$ is completely determined by the action of $H$. In particular,  for each $\chi \in \widehat{N}$, we have  $G_\chi = N \rtimes H_{\chi}$.

For each $\chi \in \widehat{N}$ and $\pi \in \widehat{H_\chi}$, let $\chi \times \pi$ denote the unitary representation of $N \rtimes H_\chi$ defined by
\begin{align*}
    (\chi \times \pi)(n, h) := \chi(n)\pi(h), \hspace{2mm}n \in N, \, h \in H_\chi.
\end{align*}

The following theorem, which is a consequence of Mackey's theory of induced representations, provides a description of the unitary dual of $G$. 

\begin{thm}[see {\cite[Theorem 4.29]{Kan-Tay}}]
\label{thm:Mackey}
    Let $G = N \rtimes H$ be a locally compact group where $N$ is an abelian locally compact group and $H$ is a locally compact group acting on $N$ such that $N$ is Mackey compatible in $G$. Let $X$ be a cross-section of the $H$-orbits in $\widehat{N}$. Then\begin{align*}
        \widehat{G} = \{\text{Ind}_{N \rtimes H_\chi}^{G}(\chi \times \pi): \pi \in \widehat{H}_{\chi}, \hspace{1mm} \chi \in X\}.
    \end{align*}
\end{thm}

For the definition of \emph{Mackey compatibility} and further details on Theorem~\ref{thm:Mackey}, we refer the reader to \cite[Section 4.3]{Kan-Tay}.

\subsection*{Type I Groups}

Let $G$ be a locally compact group. A unitary representation $\pi$ of $G$ is called \emph{multiplicity-free} if for each nontrivial decomposition $\pi = \pi_1 \oplus \pi_2 $, the unitary representations $\pi_1$ and $\pi_2$ are disjoint. The representation $\pi$ of $G$ is called \emph{type I} if $\pi$ is quasi-equivalent to a multiplicity-free representation. We say $G$ is a \emph{type I} group when every unitary representation of $G$ is type I.

\begin{prop}[see {\cite[Proposition 6.D.5]{Bekka}}]
\label{prop:Bekka}
    Let $G$ be a second countable locally compact group of type I and $\pi$ a multiplicity-free representation of $G$ in a separable Hilbert space. Then there exists a positive Borel measure $\nu$ in $\widehat{G}$ such that $\pi$ is equivalent to the direct integral $\int^{\oplus}_{\widehat{G}}\sigma d\nu(\sigma)$ and the measure $\nu$ is unique up to equivalence. 
\end{prop}

\begin{remark}
    Let $G$ be a second countable locally compact group of type I and $\rho$ a cyclic unitary representation of $G$. It follows from \cite[Chapter VII, Lemma 8.6]{Koo-RLC} that $\rho$ admits a quasi-equivalent multiplicity-free subrepresentation $\widetilde{\rho}$. Since $G$ is second countable and $\rho$ is cyclic, the Hilbert space of $\rho$ is separable. Then, the Hilbert space of $\widetilde{\rho}$ is also separable. Hence, $\widetilde{\rho}$ satisfies the assumptions of Proposition~\ref{prop:Bekka}.
\end{remark}

\section{Coincidence of the BG and KLQ \texorpdfstring{$L^p$}{Lp}-Fourier--Stieltjes Algebras}

\subsection*{Kunze--Stein Groups}

In the follwoing theorem, we establish the coincidence of the BG and KLQ $L^p$-Fourier--Stieltjes algebras for a class of nonamenable locally compact groups with the Kunze--Stein property. Since we do not use the Kunze--Stein property directly, we refer the reader to \cite[Subsection 5.1]{SW-Ex-Fnc} for its definition and further details.

\begin{thm} \label{thm: 3.2}
    Let $G$ be a nonamenable locally compact group with the Kunze--Stein property possessing a conditionally negative definite continuous function $\kappa : G \to [0, \infty)$. For each $t > 0$, let $\phi_t: G \to (0, 1]$ be defined by $\phi_t(g) := e^{-t \kappa(g)}$. We assume that
    \begin{enumerate}[(\roman*)]
        \item for each $t > 0$, there exists $p \in [1,\infty)$ 
        such that $\phi_t \in L^p(G)$;

        \item every nontrivial proper weak$^*$-closed 
        $G$-invariant ideal of $B(G)$ is of the form 
        $B_{L^q,\text{BG}}(G)$ for some $q \in [2,\infty)$.
    \end{enumerate}
    Then, for each $p \in [2,\infty)$, we have
    \begin{align*}
        B_{L^p,\text{BG}}(G) = B_{L^p,\text{KLQ}}(G).
    \end{align*}
\end{thm}

\begin{proof}
First, we recall that $B_{L^2, \text{KLQ}}(G) = B_r(G) =  B_{L^2, \text{BG}}(G)$ (see \cite[Lemma 2.10]{FTW}). Thus, it suffices to consider the case $p \in (2, \infty)$. Let $t_0 := \text{inf}\{t > 0 : \phi_t \in L^1(G)\}$. Since $G$ is nonamenable, we have $t_0 > 0$ (see the proof of \cite[Theorem 3.3]{FTW}). It follows from the proof of \cite[Corollary 5.7]{SW-Ex-Fnc} that, for each $q \in [2, \infty)$, $\phi_t \in B_{L^q, \text{BG}}(G) \Longleftrightarrow t \geq \frac{t_0}{q}$. Now, let $p \in (2, \infty)$ be fixed. Since $B_{L^p, \text{BG}}(G) \subset B_{L^p, \text{KLQ}}(G)$, for all $t \geq \frac{t_0}{p}$, $\phi_t \in B_{L^p, \text{KLQ}}(G)$ . On the other hand, it follows from the proof of \cite[Theorem 3.3]{FTW} for the Kunze--Stein case that $\phi_t \notin B_{L^p, \text{KLQ}}(G)$ when $t < \frac{t_0}{p}$. Hence, $\phi_t \in B_{L^p, \text{KLQ}}(G) \Longleftrightarrow t \geq \frac{t_0}{p}$.

In particular, $B_{L^p, \text{KLQ}}(G)$ is a nontrivial proper weak$^*$-closed $G$-invariant ideal of $B(G)$. Then, by assumption (ii), there exists $q \in [p, \infty)$ such that $B_{L^p, \text{KLQ}}(G) = B_{L^q, \text{BG}}(G)$. By the above characterization, we have $t \geq \frac{t_0}{p} \Longleftrightarrow t \geq \frac{t_0}{q},$ for every $t > 0$. Since $t_0 > 0$, it follows that $p = q$. Hence,
$B_{L^p, \text{BG}}(G) = B_{L^p, \text{KLQ}}(G),$ for each $p \in [2, \infty)$.

\end{proof}

\begin{remark}
The assumptions of the preceding theorem are satisfied by the groups $\text{SL}(2, \mathcal{K})$, where $\mathcal{K} = \mathbb{R}, \mathbb{C}$, or $\mathbb{Q}_p$, and noncompact closed subgroups of $\text{Aut}(T)$ acting transitively on the vertices and boundary of a homogeneous tree $T$ satisfying Tits's independence property. 

These groups have the Kunze--Stein property (see \cite[Subsection 5.1]{SW-Ex-Fnc}) and satisfy the integrability condition in assumption (i) (see \cite[Section 4]{SW-Ex-Fnc}). The characterization of weak$^*$-closed $G$-invariant ideals required in assumption (ii) follows from \cite[Theorem 7.3]{WRS-Fnc} for $\text{SL}(2, \mathbb{R})$, \cite[Theorem A]{de Laat} for $\text{SL}(2, \mathbb{C})$, and \cite[Theorems 4.24 and A]{HLS-Tr-IMRN} for $\text{SL}(2, \mathbb{Q}_p)$ and the mentioned subgroups of Aut$(T)$, respectively. Moreover, \cite[Proposition 4.2]{SW-Ex-Fnc} allows us to pass from $L^{p+}$-Fourier--Stieltjes algebras to $L^p$-Fourier--Stieltjes algebras whenever required.

In particular, Theorem~\ref{thm: 3.2} recovers the coincidence result for $\text{SL}(2, \mathbb{R})$ established in \cite[Theorem 7.1]{FTW}.
\end{remark}

\subsection*{Semidirect Product Groups}

In the remainder of this section, we investigate the BG and KLQ $L^p$-Fourier--Stieltjes algebras for a class of semidirect product groups by identifying their unitary duals.

Let $G = N \rtimes H$ be a second countable locally compact group where $N$ is a noncompact abelian locally compact group, which is regularly embedded in $G$, and $H$ is a nonamenable type I group. Then $N$ is Mackey compatible (see \cite[Section 4.3]{Kan-Tay} for details). Let $X$ be a cross-section of the $H$-orbits in $\widehat{N}$. We define $X_0 := \{\chi \in X: H_\chi = H\}$ and assume that for every $\chi \in X \setminus X_0$, the stabilizer group $H_\chi$ is an amenable type I group. From now on, we assume that $G = N \rtimes H$ satisfies the above assumptions, unless otherwise stated.

\begin{prop}
    Let $G = N \rtimes H$. Then $G$ is type I.
\end{prop}

\begin{proof}
    By assumption, the stabilizer group $H_\chi$ is type I for each $\chi \in X$. We also note that $G_\chi = N \rtimes H_\chi$ and $\chi': G_\chi \to \mathbb{T}$ defined by $\chi'(n, h) = \chi(n)$ for all $(n, h) \in G_\chi$ is a unitary representation since $H_\chi$ stabilizes $\chi$. Then, it follows from \cite[Proposition 2.36]{Carter} that  $G$ is a type I group.
\end{proof}

\begin{prop} \label{prop:3.10}
    Let $G = N \rtimes H$. Then 
    \begin{align*}
          \widehat{G} = \bigcup_{\chi \in X_0}\{\chi \times \pi: \pi \in \widehat{H} \}\cup \{\text{Ind}_{N \rtimes H_\chi}^G(\chi \times \pi): \pi \in \widehat{H_\chi}, \, \chi \in X \setminus X_0\}.
    \end{align*}
\end{prop}

\begin{proof}
    The result follows from Theorem~\ref{thm:Mackey}.
\end{proof}

In the following, we use the notation $\pi \prec \sigma$ to denote weak containment of unitary representations.

\begin{prop} \label{prop:3.11}
    Let $G = N \rtimes H$. Then: 
    \begin{align*}
     \widehat{G}_r = \bigcup_{\chi \in X_0}\{\chi \times \pi: \pi \in \widehat{H}_r \}\cup \{\text{Ind}_{N \rtimes H_\chi}^G(\chi \times \pi): \pi \in \widehat{H_\chi}, \, \chi \in X \setminus X_0\}
    \end{align*}
    where $\widehat{G}_r$ and $\widehat{H}_r$ denote the reduced unitary duals of $G$ and $H$, respectively.
    
\end{prop}

\begin{proof}
    Let $\chi \in X\setminus X_0$ and $K_\chi := N \rtimes H_\chi$. Then, $K_\chi$ is amenable which implies that 
    $\chi \times\pi \prec \lambda_{K_\chi}$ for each $\pi \in \widehat{H_\chi}$. Hence, $\text{Ind}_{K_\chi}^G (\chi \times \pi) \prec \lambda_G$ for each $\pi \in \widehat{H_\chi}$.
    
     Assume that for some $\chi \in X_0$ and $\pi \in \widehat{H}$, we have $(\chi \times \pi) \prec \lambda_G$, then $\text{Res}_H^G(\chi \times \pi) \prec \lambda_H$ which implies that $\pi \prec \lambda_H$ and consequently $\pi \in \widehat{H}_r$. Conversely, assume that $\pi \in \widehat{H}_r$. We note that for each $\chi \in X_0$, $\chi \prec \lambda_N$, and the induced representation $\text{Ind}_N^G(\chi)$ is unitarily equivalent to $ (\chi \times \lambda_H)$. Thus, 
    $(\chi \times \pi) \prec \lambda_G$,
    since $\pi \prec \lambda_H$.
\end{proof}

\begin{remark}
Let $\mathcal{K}$ be a local field. Then $\text{SL}(n,\mathcal{K})$ is a type I group for every $n \geq 2$ (see \cite{Bek-Ech}).
\end{remark}

\begin{eg} \label{example: 3.13}
Let $G = \mathcal{K}^2 \rtimes \text{SL}(2, \mathcal{K})$, where $\mathcal{K}$ is a local field and $\text{SL}(2, \mathcal{K})$ acts on $N = \mathcal{K}^2$ by the natural linear action. Fix a nontrivial unitary character $\psi$ of $\mathcal{K}$. Considering that the additive group $\mathcal{K}$ is self-dual, we identify $\widehat{N}$ with $\mathcal{K}^2$ via characters $\chi_{\xi}(n) = \psi(\langle n, \xi \rangle), \hspace{2mm} \xi \, , n \in \mathcal{K}^2$. Then the action of $H = \text{SL}(2, \mathcal{K})$ on $\widehat{N}$ is given by $h\cdot \xi = h^{-T}\xi.$ Since $\text{SL}(2, \mathcal{K})$ acts transitively on $\mathcal{K}^2 \setminus \{0\}$, we have that $\widehat{N}/H = \{\{0\}, \mathcal{K}^2 \setminus\{0\}\}.$ Let $X = \{0, e_1\}$ be a cross-section of the $H$-orbits in $\widehat{N}$, where $e_1 = (1, 0)^T$. Then, $H_0 = \text{SL}(2,\mathcal{K})$ and $H_{e_1} = \left\{\begin{pmatrix}
    1 & 0 \\ k & 1
\end{pmatrix}: k \in \mathcal{K}\right\},$ which is isomorphic to $\mathcal{K}$ as a topological group. By Theorem~\ref{thm:Mackey}, we have 
\begin{align*}
    \widehat{G} = \widehat{\text{SL}(2, \mathcal{K})} \circ q \hspace{1mm}\cup \hspace{1mm} \{\text{Ind}_{\mathcal{K}^2 \rtimes H_{e_1}}^G (\chi_{e_1} \times \phi): \phi \in \widehat{H_{e_1}}\}
\end{align*}
where $q: G \to H$ is the canonical quotient map.

Let $L := \mathcal{K}^2 \rtimes H_{e_1}$. Then $L$ is an amenable group which implies that $\pi \prec \lambda_L$, for each unitary representation $\pi$ of $L$. Therefore, $(\chi_{e_1} \times \phi) \prec \lambda_L$ for each $\phi \in \widehat{H_{e_1}}$ which implies that $\text{Ind}_{\mathcal{K}^2 \rtimes H_{e_1}}^G (\chi_{e_1} \times \phi) \prec \lambda_G$ for each $\phi \in \widehat{H_{e_1}}$.

Since $N$ is amenable, we have $1_N \prec \lambda_N$ and this implies that $(\text{Ind}_N^G 1_N) \prec (\text{Ind}_N^G \lambda_N)$.
 We also note that $\lambda_G$ is unitarily equivalent to $(\text{Ind}_N^G \lambda_N)$ and $(\text{Ind}_N^G 1_N)$  
 is unitarily equivalent to $\lambda_H \circ q$. Hence, we have that $\lambda_H \circ q \prec \lambda_G$. Conversely, if $\pi \circ q \prec \lambda_G$, then by restricting to $H$, we have $\pi \prec \lambda_H$. Therefore,
\begin{align*}
\widehat{G}\setminus\widehat{G}_r = \{\pi \circ q : \pi \in \widehat{\text{SL}(2, \mathcal{K})} \setminus \widehat{\text{SL}(2, \mathcal{K})_r}\}.
\end{align*}

\end{eg}

\begin{eg} \label{eg:3.14}
Let $H := \mathcal{K}^2 \rtimes \text{SL}(2,\mathcal{K})$, and fix a nontrivial unitary character $\psi$ of $\mathcal K$. We realize $H$ as the closed subgroup of $\text{SL}(3, \mathcal{K})$ given by 
\begin{align*}
    (u, A) \hookrightarrow \begin{pmatrix}
        A & u \\ 0 & 1
    \end{pmatrix},
\end{align*}
where $u \in \mathcal{K}^2$ and $A \in \text{SL}(2, \mathcal{K})$. We identify $N := \mathcal{K}^3$ with $\mathcal{K}^2 \oplus \mathcal{K}$, and write its elements as $(\xi, k)$, where $\xi \in \mathcal{K}^2$ and $k \in \mathcal{K}$. The action of $H$ on $N$ is given by the ordinary matrix multiplication $(u, A)\cdot (\xi, k) = (A\xi + ku, \,k).$ This defines the semidirect product group $G = N \rtimes H$. We identify the unitary dual $\widehat{N}$ with $N$ via the correspondence $(\xi, k) \in \mathcal{K}^2 \oplus \mathcal{K} \longleftrightarrow \chi_{(\xi, k)} \in \widehat{N},$ where $\chi_{(\xi, k)}(\eta, t) = \psi(\langle \eta, \xi \rangle + tk), \hspace{1mm}\eta \in \mathcal{K}^2, \, t \in \mathcal{K}$. Then, the action of $H$ on $\widehat{N}$ is given by 
\begin{align*}
    (u, A) \cdot (\xi, k) = (A^{-T}\xi, k - \langle u , A^{-T}\xi \rangle).
\end{align*}
If $\xi = 0$, then $(u, A) \cdot (0, k) = (0, k)$ for every $(u, A) \in H$. Thus each point $(0, k)$, $k \in \mathcal{K}$, forms a singleton orbit. Now let $\xi \neq 0$. Since $\text{SL}(2, \mathcal{K})$ acts transitively on $\mathcal{K}^2 \setminus \{0\}$, we may choose $A$ in $\text{SL}(2, \mathcal{K})$ such that $A^{-T}\xi = e_1$, where $e_1 = (1, 0)^T$. Moreover, by choosing $u \in \mathcal{K}^2$ appropriately, the second coordinate can be shifted arbitrarily. It follows that every point $(\xi, k)$ with $\xi \neq 0$ lies in the $H$-orbit of $(e_1, 0)$. Therefore, $\widehat{N}/H = \{(0, k): k \in \mathcal{K}\} \cup \{H(e_1, 0)\},$
where $H(e_1, 0) = (\mathcal{K}^2 \setminus \{0\}) \times \mathcal{K}$. Hence, a cross-section for the $H$-orbits is given by $X = \{(0, k): k \in \mathcal{K}\} \cup \{(e_1, 0)\}.$

For each $k \in \mathcal{K}$, $H_{(0, k)} = H$ and 
\begin{align*}
    H_{(e_1, 0)} = \{(u, A) \in H: A^{-T}e_1 = e_1, \langle u, e_1 \rangle = 0\} = \left\{\begin{pmatrix}
        1 & 0 & 0 \\ t & 1 & s\\ 0 & 0& 1
    \end{pmatrix}: t, s \in \mathcal{K}\right\},
\end{align*}
which is isomorphic to $\mathcal{K}^2$ as a topological group. Therefore, Theorem~\ref{thm:Mackey} gives
\begin{align*}
    \widehat{G} = \{\chi_{(0, k)} \times \pi: k \in \mathcal{K}, \, \pi \in \widehat{H}\} \cup \{\text{Ind}_{N \rtimes H_{(e_1, 0)}}^G (\chi_{(e_1, 0)} \times \tau): \tau \in \widehat{H_{(e_1, 0)}}\}.
\end{align*}
Since the stabilizer subgroup $H_{(e_1, 0)}$ is abelian and hence amenable, every irreducible representation $\tau \in \widehat{H_{(e_1, 0)}}$ is weakly contained in the left regular representation of $H_{(e_1, 0)}$. It follows that 
\begin{align*}
    \text{Ind}_{N \rtimes H_{(e_1, 0)}}^G(\chi_{(e_1, 0)} \times \tau ) \prec \lambda_G.
\end{align*}
Therefore, the irreducible representations of $G$ that are not weakly contained in $\lambda_G$ must come from the fixed characters $(0, k)$, $k \in \mathcal{K}$. Consequently,
\begin{align*}
    \widehat{G} \setminus \widehat{G}_r = \bigcup _{k \in \mathcal{K}}\{\chi_{(0, k)}\times \pi : \pi \in \widehat{H} \setminus \widehat{H}_r\}.
\end{align*}
In particular, $\widehat{G} \setminus \widehat{G}_r$ consists of a family of copies of $\widehat{H} \setminus \widehat{H}_r$, indexed by the fixed characters of the action $H$ on $\widehat{N}$.
\end{eg}

Although the above examples both satisfy the assumptions discussed so far, their orbit structures in $\widehat{N}$ have  different behavior.

\begin{eg}
Let $G = \mathcal{K}^3 \rtimes \text{SL}(3, \mathcal{K})$, where $\mathcal{K}$ is a local field and $\text{SL}(3,\mathcal{K})$ acts on $N = \mathcal{K}^3$ by the natural linear action. Fix a nontrivial unitary character $\psi$ of $\mathcal{K}$. We identify $\widehat{N}$ with $\mathcal{K}^3$ via characters
$\chi_\xi(n)=\psi(\langle n,\xi \rangle), \hspace{1mm} \xi\in \mathcal{K}^3,$ $n\in N$. The induced action of $H = \text{SL}(3,\mathcal{K})$ on $\widehat{N}$ is given by
$(A \cdot \chi)(n) = \chi (A^{-1}n), \hspace{1mm} A \in \text{SL}(3,\mathcal{K}), \ \chi \in \widehat{N}, \ n \in N$. Under the identification of $\widehat{N}$ with $\mathcal{K}^3$, this corresponds to $A \cdot \xi = A^{-T}\xi$. Since $\text{SL}(3,\mathcal{K})$ acts transitively on $\mathcal{K}^3\setminus\{0\}$, we have $\widehat{N}/H = \{\{0\}, \mathcal{K}^3\setminus\{0\}\}$. Let $X = \{0,e_1\}$ be a cross-section of the $H$-orbits in $\widehat{N}$, where $e_1 = (1,0,0)^T$. Then $H_0 = \text{SL}(3,\mathcal{K})$, and $H_{e_1} = \{A\in \text{SL}(3,\mathcal{K}):A^{-T}e_1 = e_1\}$.
Equivalently,
\begin{align*}
H_{e_1} = \left\{
\begin{pmatrix}
1 & 0 & 0\\
k_1 & a & b\\
k_2 & c & d
\end{pmatrix}: \begin{pmatrix}
a & b\\
c & d
\end{pmatrix}\in \text{SL}(2,\mathcal{K}),\ k_1,k_2\in \mathcal{K}
\right\},
\end{align*}
which is isomorphic to $\mathcal{K}^2\rtimes \text{SL}(2,\mathcal{K})$ as a topological group. By Theorem~\ref{thm:Mackey}, we obtain
\begin{align*}
\widehat{G} = \widehat{\text{SL}(3,\mathcal{K})}\circ q \cup
\left\{
\text{Ind}^{G}_{\mathcal{K}^3\rtimes H_{e_1}}
(\chi_{e_1}\times \pi) :
\pi\in \widehat{H}_{e_1}
\right\},
\end{align*}
where $q:G\to \text{SL}(3,\mathcal{K})$ is the canonical quotient map. We note that $H_{e_1}$ is nonamenable, since it contains $\text{SL}(2,\mathcal{K})$ as a closed subgroup. Thus, unlike the groups considered so far, the stabilizer subgroup corresponding to the nontrivial orbit is not amenable. Hence, the argument of Proposition~\ref{prop:3.11} does not apply in this case. In particular, we can no longer conclude that
\begin{align*}
\text{Ind}^{G}_{\mathcal{K}^3\rtimes H_{e_1}}
(\chi_{e_1}\times \pi)
\prec \lambda_G,
\end{align*}
for every $\pi\in \widehat{H}_{e_1}$. Hence, the nontrivial orbit may contribute irreducible representations that are not weakly contained in the left regular representation of $G$. This shows that the assumption requiring the stabilizer groups $H_\chi$, $\chi\in X\setminus X_0$, to be amenable is essential in our results.
\end{eg}
  
\begin{prop}\label{thm:3.11}
    Let $G = N \rtimes H$ and suppose that $X_0$ is countable. Then 
    \begin{align*}
         B_r(G) = \overline{B_0(G)}^{\sigma(B(G), C^*(G))}.
    \end{align*}
\end{prop}

\begin{proof}
    Let $\rho$ be a cyclic unitary representation of $G$ such that $A_{\rho}(G) \subset B_0(G)$. Since $G$ is a type I group, $\rho$ admits a quasi-equivalent multiplicity-free subrepresentation $\widetilde{\rho}$. Then, $A_{\rho}(G) = A_{\widetilde{\rho}}(G)$. By Proposition~\ref{prop:Bekka}, there exists a positive Borel measure $\nu$ on $\widehat{G}$ such that $\widetilde{\rho}$ is equivalent to the direct integral representation $\int^{\oplus}_{\widehat{G}}\sigma d\nu(\sigma)$.
    
    Fix $\chi \in X_0$. Since $G$ and $H$ are second countable type I groups, their unitary duals $\widehat{G}$ and $\widehat{H}$, equipped with their canonical Borel structures, are standard Borel spaces. We define the map $\phi_\chi: \widehat{H} \to \widehat{G}$ via $\phi_\chi(\pi) = \chi \times \pi$. Then $\phi_\chi$ is a continuous injection and hence is a Borel map. It follows from \cite[Corollary 15.2]{Kechris} that $\phi_\chi(\widehat{H}) = \{\chi \times \pi: \pi \in \widehat{H}\}= \chi \times \widehat{H} $ is a Borel subset of $\widehat{G}$.
    
    If we assume that $\nu (\chi \times \widehat{H}) > 0$ for some $\chi \in X_0$, then the representation $\widetilde{\rho}_\chi := \int_{\chi \times \widehat{H}}^{\oplus}\sigma d \nu(\sigma)$ is a nonzero subrepresentation of $\widetilde{\rho}$ which implies that $A_{\widetilde{\rho}_\chi}(G) \subset A_{\widetilde{\rho}}(G)$. We also note that $A_{\widetilde{\rho}}(G) \subset B_0(G)$ since $A_{\widetilde{\rho}}(G) = A_{\rho}(G)$. On the other hand, representations of the form $\chi \times \pi$, where $\pi$ is a unitary representation of $H$ and $\chi \in X_0$, give rise to coefficient functions of the form
    \begin{align*}
        (\chi \times \pi)_{\xi, \eta}(n, h) = \chi(n) 
\langle \pi(h)\xi, \eta \rangle, \text{ for } \xi, \eta \in \mathcal{H}_\pi.
    \end{align*}
    We claim that every nonzero coefficient function of this form does not belong to $C_0(G)$. Indeed, suppose that $f$ is a nonzero coefficient function of this form. Then there exists $h_0 \in H$ such that $f(n, h_0) \neq 0$ for all $n \in N$.
    Let $K \subset G$ be any compact set. Since $N$ is noncompact, there exists $n_0 \in N$ such that $(n_0, h_0) \notin K$. It follows that 
    $|f(n_0, h_0)| > 0$. Thus, $f \notin C_0(G)$, and then $f \notin B_0(G)$. This contradicts the fact that $A_{\widetilde{\rho}_\chi}(G) \subset B_0(G)$. Hence, $\nu (\chi \times \widehat{H}) = 0$, for every $\chi \in X_0$. Since $X_0$ is countable, we have $\nu (\bigcup_{\chi \in X_0} (\chi \times \widehat{H})) \leq \sum_{\chi \in X_0} \nu(\chi \times \widehat{H}) = 0.$ Also, by Proposition~\ref{prop:3.11}, $\widehat{G} \setminus \widehat{G}_r \subset \bigcup_{\chi \in X_0} (\chi \times \widehat{H})$. Hence, $\nu (\widehat{G} \setminus \widehat{G}_r) = 0$, and almost every irreducible representation in the direct integral decomposition of $\widetilde{\rho}$ is weakly contained in $\lambda_G$. Hence, $B_{\rho}(G) = B_{\widetilde{\rho}}(G) \subset B_r(G)$. Since $\rho$ was an arbitrary cyclic representation of $G$ with $A_{\rho}(G) \subset B_0(G)$, $B_r(G) = \overline{A(G)}^{\sigma(B(G), C^*(G))}$, and $A(G) \subset B_0(G)$, we can conclude that
    \begin{align*}
        B_r(G) = \overline{B_0(G)}^{\sigma(B(G), C^*(G))}.
    \end{align*}    
\end{proof}

\begin{prop}\label{prop:2}
    Let $G = N \rtimes H$ and suppose that $B_0(G)|_H \subset B_r(H)$. Then
    \begin{align*}              \overline{B_0(G)}^{\sigma(B(G) , C^*(G))} = B_r(G).
    \end{align*}
\end{prop}

\begin{proof}
    Let $\rho$ be a cyclic unitary representation of $G$ such that $A_\rho (G) \subset B_0(G)$. Since $H$ is a closed subgroup of $G$, we have
    \begin{align*}
        A_{\text{Res}_H^G (\rho)}(H) = A_\rho (G)|_H \subset B_0(G)|_H  \subset B_r(H).
    \end{align*}
    Hence, $\text{Res}_H ^G(\rho) \prec \lambda_H$. Since $G$ is a second countable type I group, $\rho$ admits a quasi-equivalent multiplicity-free subrepresentation $\widetilde{\rho}$. By Proposition~\ref{prop:Bekka}, there exists a positive Borel measure $\nu$ on $\widehat{G}$ such that $\widetilde{\rho}$ is equivalent to the direct integral representation $\int^{\oplus}_{\widehat{G}}\sigma d\nu(\sigma)$. Since $\rho$ and $\widetilde{\rho}$ are quasi-equivalent, we also have $\text{Res}_H^G(\widetilde{\rho}) \prec \lambda_H$. Since restriction to a closed subgroup commutes with direct integrals, we have
    \begin{align*}
    \text{Res}_H^G(\widetilde{\rho}) \cong \int_{\widehat{G}}^{\oplus} \text{Res}_H^G (\sigma) \hspace{1mm} d\nu(\sigma). 
    \end{align*}
    
     Since $\text{Res}_H^G (\widetilde{\rho}) \prec \lambda_H$, we have $\text{Res}_H^G (\sigma) \prec \lambda_H$, for $\nu$-almost every $\sigma \in \widehat{G}$. We now assume that $\sigma \in \widehat{G} \setminus \widehat{G}_r$. By Proposition~\ref{prop:3.11}, $\sigma = \chi \times \pi$ for some $\chi \in X_0$ and $\pi \in \widehat{H} \setminus \widehat{H}_r$. Then, $\text{Res}_H^G(\sigma) = \pi$ which is not weakly contained in $\lambda_H$. This implies that $\nu (\widehat{G} \setminus \widehat{G}_r) = 0$. Then, $\nu$-almost every representation in the direct integral decomposition of $\widetilde{\rho}$ is weakly contained in $\lambda_G$.  Hence, $\widetilde{\rho} \prec \lambda_G$ and since $\rho$ and $\widetilde{\rho}$ are quasi-equivalent, $\rho \prec \lambda_G$. This implies $B_{\rho}(G) \subset B_r(G)$. Since $\rho$ was an arbitrary cyclic representation of $G$ with $A_{\rho}(G) \subset B_0(G)$, we can conclude that $B_r(G)$ coincides with $\overline{B_0(G)}^{\sigma(B(G), C^*(G))} $.
\end{proof}

\begin{thm}\label{main theorem}
    Let $G = N \rtimes H$. Suppose that at least one of the following conditions holds:
    \begin{enumerate}[(\roman*)]
        \item $X_0$ is countable.
        \item $B_0(G)|_H \subset B_r(H)$.
    \end{enumerate}  
    Then, for each $p \in [2, \infty)$ we have 
        \begin{align*}
            B_{L^p, \text{BG}}(G) = B_r(G) = B_{L^p, \text{KLQ}}(G).
        \end{align*}
\end{thm}

\begin{proof}
     For each $p \in [2, \infty)$,  it is known that the Fourier algebras $A_{L^p,\text{BG}}(G)$ and $A_{L^p,\text{KLQ}}(G)$ are contained in the Rajchman algebra $B_0(G)$. Then, Propositions~\ref{thm:3.11} and~\ref{prop:2} with the monotonicity of $L^p$-Fourier--Stieltjes algebras with respect to inclusion imply that 
     \begin{align*}
     B_{L^p, \text{BG}}(G) = B_r(G) = B_{L^p, \text{KLQ}}(G), \hspace{1mm}\text{for each} \hspace{2mm}p \in [1, \infty).
     \end{align*}
\end{proof}

\begin{remark} 
    \begin{enumerate}[(\roman*)]
        \item Let $G = \mathcal{K}^2 \rtimes \text{SL}(2,\mathcal{K}).$ By Example~\ref{example: 3.13}, the set $X_0$ consists only of the
        trivial character and is therefore countable. Hence, Proposition~\ref{thm:3.11} implies that
        \begin{align*}
            \overline{B_0(G)}^{\sigma(B(G),C^*(G))} = B_r(G).
        \end{align*}

        \item Let $G = \mathcal{K}^3 \rtimes H$ where $H = \mathcal{K}^2 \rtimes \text{SL}(2,\mathcal{K})$.
        By Example~\ref{eg:3.14}, the set $X_0 = \{\chi_{(0,k)}: k \in \mathcal{K}\}$ is uncountable, so Proposition~\ref{thm:3.11} is not applicable. However, by part (i), we have
        \begin{align*}
            \overline{B_0(H)}^{\sigma(B(H),C^*(H))} = B_r(H).
        \end{align*}
        Hence, $B_0(G)|_H \subset B_0(H)\subset B_r(H)$, and Proposition~\ref{prop:2} implies that
        \begin{align*}
            \overline{B_0(G)}^{\sigma(B(G),C^*(G))} = B_r(G).
        \end{align*}
    \end{enumerate}

    In particular, the conclusion of Theorem~\ref{main theorem} holds for both of these groups.
\end{remark}

\section{Closed Ideals of Unitary Dual for Semidirect Product Groups}

We recall that, for a locally compact group $G$, a subset $S \subset \widehat{G}$ is called an \emph{ideal} of $\widehat{G}$ if, for every $\pi \in S$ and $\rho \in \widehat{G}$, the unitary representation $\pi \otimes \rho$ is weakly contained in $S$. An ideal $S$ is said to be \emph{closed} if it is closed in the Fell topology on $\widehat{G}$.

In the following, we describe the closed ideals of the unitary dual $\widehat{G}$ in terms of those of $\widehat{H}$, where $G = N \rtimes H$ satisfies the assumptions stated in the previous section.

For a closed ideal $S \subset \widehat{H}$, we define
\begin{align*}
    S_G := \widehat{G}_r \cup  \bigcup_{\chi \in X_0}\{\chi \times \pi: \pi \in S \setminus \widehat{H}_r\} \subset \widehat{G}.
\end{align*}

\begin{thm} \label{thm:3.19}
    Let $G = N \rtimes H$.  Then for every closed ideal $S \subset \widehat{H}$, the set $S_G$ is a closed ideal of $\widehat{G}$. Moreover, every closed ideal $J \subset \widehat{G}$ is of this form for a unique closed ideal $S \subset \widehat{H}$.
\end{thm}

\begin{proof}
    Let $S \subset \widehat{H}$ be a closed ideal. We first show that $S_G$ is an ideal of $\widehat{G}$. Let $\gamma \in S_G$ and $\rho \in \widehat{G}$. If $\gamma \in \widehat{G}_r$, then $\gamma \prec \lambda_G$, and by Fell's absorption principle, $\gamma \otimes \rho \prec \lambda_G$. Hence, $\gamma \otimes \rho \prec S_G$. Now we assume that $\gamma = \chi_1 \times \pi$, where $\chi_1 \in X_0$ and $\pi \in S \setminus \widehat{H}_r$. If $\rho \prec \lambda_G$, then again by Fell's absorption principle, $(\chi_1 \times \pi) \otimes \rho \prec \lambda_G$, and so $(\chi_1 \times \pi) \otimes \rho \prec S_G$. Otherwise, by Proposition~\ref{prop:3.10}, $\rho = \chi_2 \times \sigma$ for some $\chi_2 \in X_0$ and $\sigma \in \widehat{H} \setminus \widehat{H}_r$. Then
    \begin{align*}
        (\chi_1 \times \pi) \otimes (\chi_2 \times \sigma) = (\chi_1 \chi_2) \times (\pi \otimes \sigma).
    \end{align*}
    Since $\chi_1, \chi_2 \in X_0$, $\chi_1 \chi_2$ also belongs to $X_0$. We also note that $\chi \times S \subset S_G$ for every $\chi \in X_0$. Since $S$ is an ideal of $\widehat{H}$, we have $\pi \otimes \sigma \prec S$ and then $(\chi_1 \chi_2) \times (\pi \otimes \sigma) \prec (\chi_1 \chi_2) \times S \prec S_G.$ Hence, $S_G$ is an ideal of $\widehat{G}$.
    
    We now show that $S_G$ is closed. Let $(\gamma_i)$ be a net in $S_G$ such that $\gamma_i \rightarrow \gamma$ in $\widehat{G}$. If eventually $\gamma_i \in \widehat{G}_r$, then $\gamma \in \widehat{G}_r \subset S_G$ since $\widehat{G}_r$ is closed. Otherwise, after passing to a subnet, we may assume that $\gamma_i = \chi_i \times \pi_i$, where $\chi_i \in X_0$ and $\pi_i \in S \setminus \widehat{H}_r$. If $\gamma \notin \widehat{G}_r$, then by Proposition~\ref{prop:3.10}, $\gamma = \chi \times \pi$ for some $\chi \in X_0$ and $\pi \in \widehat{H} \setminus \widehat{H}_r$. Since $\gamma_i \to \gamma$, for every subnet $(\gamma_j)$ of $(\gamma_i)$, we have $\gamma \prec \bigoplus_{j} \gamma_j$. Then $\pi = \text{Res}_H ^G (\gamma) \prec \bigoplus_j \text{Res}_H ^G(\gamma_j) = \bigoplus_j \pi_j$ and this implies that $\pi \in \overline{S} = S$. Hence, $\gamma \in S_G$ and $S_G$ is closed.
    
    Let $J \subset \widehat{G}$ be a closed ideal. We note that the trivial character $1_N$ belongs to $X_0$. Now, let $S := \{\pi\in\widehat{H} : 1_N \times \pi \in J\}$. We first show that $\chi \times \pi \in J$ if and only if $\pi \in S$, for every $\chi \in X_0$ and $\pi \in \widehat{H}$. Let $\chi\in X_0$ and $\pi \in \widehat{H}$. Since $\chi$ is fixed by $H$, its inverse $\chi^{-1}$ is also fixed by $H$, and then $\chi^{-1}\in X_0$. If $\chi\times\pi\in J$, then, since $J$ is an ideal, $(\chi \times \pi) \otimes (\chi^{-1} \times 1_H) = 1_N \times \pi \in J$. Conversely, if $1_N \times \pi \in J$, then $(1_N \times \pi) \otimes (\chi \times 1_H) = \chi \times \pi \in J$. Hence, for every $\chi \in X_0$ and $\pi \in \widehat{H}$, $\chi \times \pi \in J$ if and only if $1_N \times \pi \in J$. Now we define
    \begin{align*}
    \phi: \widehat{H}\to \widehat{G},
    \hspace{3mm} \phi(\pi):= 1_N \times \pi.
    \end{align*}
    Then $S = \phi^{-1}(J)$. Since $\phi$ is continuous and $J$ is closed, $S$ is closed in $\widehat{H}$. 
    
    We next show that $S$ is an ideal of $\widehat{H}$. Let $\pi\in S$ and $\sigma \in \widehat{H}$. Then $1_N \times \pi \in J$. Since $1_N \times \sigma \in \widehat{G}$ and $J$ is an ideal,
    \begin{align*}
    (1_N \times \pi) \otimes (1_N \times \sigma) = 1_N \times (\pi \otimes \sigma) \prec J.
    \end{align*}
    We note that if $\gamma \in \widehat{H}$ such that $\gamma \prec \pi \otimes \sigma$, then $1_N \times \gamma \prec J$. Since $J$ is closed, $1_N \times \gamma \in J$ and then $\gamma \in S$. Thus, $\pi \otimes \sigma \prec S$, and $S$ is an ideal of $\widehat{H}$.
    
    We claim that $J = S_G$. The inclusion $S_G \subset J$ follows from the equivalence above and the fact that $\widehat{G}_r \subset J$. Conversely, let $\gamma \in J$. If $\gamma \in \widehat{G}_r$, then $\gamma \in S_G$. Otherwise, by Proposition~\ref{prop:3.11}, $\gamma = \chi \times \pi$ for some $\chi \in X_0$ and $\pi \in \widehat{H} \setminus \widehat{H}_r$. Since $\gamma \in J$, the equivalence above gives $1_N \times \pi \in J$, and then $\pi \in S$. This implies that $\gamma \in S_G$, and $J = S_G$.
    
    To prove uniqueness, we assume that $J = S_G = T_G$ for closed ideals $S, \hspace{1mm} T \subset \widehat{H}$. Let $\pi \in \widehat{H} \setminus \widehat{H}_r$. We have
    \begin{align*}
    \pi \in S \Longleftrightarrow 1_N \times \pi \in J \Longleftrightarrow \pi \in T.
    \end{align*}
    Thus, $S \setminus \widehat{H}_r = T \setminus \widehat{H}_r$. Since both $S$ and $T$ contain $\widehat{H}_r$, it follows that $S = T$.
    \end{proof}

Let $H = \text{SL}(2,\mathcal{K})$, where $\mathcal{K} = \mathbb{R}$, $\mathbb{C}$, or $\mathbb{Q}_p$. For each $p \in [2,\infty]$, let
\begin{align*}
\widehat{H}_{L^p} := \{\sigma \in \widehat{H} : \sigma \text{ is an } L^p\text{-representation of } H\}.
\end{align*}
It is known that $\widehat{H}_{L^p}$ is a closed ideal of $\widehat{H}$ for each $p \in [2,\infty]$ and that every closed ideal of $\widehat{H}$ is of this form. (see \cite{de Laat}, \cite{HLS-Tr-IMRN})

\begin{prop} \label{cor:3.20}
    Let $G = \mathcal{K}^2 \rtimes H$ and $q: G \to H$ denote the quotient map. For each $p \in [2,\infty]$, define
    \begin{align*}
        S_p := \widehat{G}_r \cup \{\sigma\circ q: \sigma \in \widehat{H}_{L^p}\setminus \widehat{H}_r\}.
    \end{align*}
    Then the closed ideals of $\widehat{G}$ are precisely the sets $S_p$, where $p\in[2,\infty]$.
\end{prop}

\begin{proof}
    By Example~\ref{example: 3.13},
    \begin{align*}
        \widehat{G} \setminus \widehat{G}_r = \{\sigma \circ q: \sigma \in \widehat{H}\setminus \widehat{H}_r\}.
    \end{align*}
    Since every closed ideal of $\widehat{H}$ is of the form $\widehat{H}_{L^p}$ for some $p\in[2,\infty]$, the result follows immediately from Theorem~\ref{thm:3.19}.
\end{proof}

\begin{prop}
    Let $G = \mathcal{K}^2 \rtimes H$. Then the weak$^*$-closed $G$-invariant ideals of $B(G)$ are precisely the Fourier--Stieltjes spaces $B_{S_p}(G)$, where $p \in [2, \infty]$.
\end{prop}

\begin{proof}
    By Proposition~\ref{cor:3.20}, the closed ideals of $\widehat{G}$ are precisely the sets $S_p$, for $p \in [2,\infty]$. Since weak$^*$-closed $G$-invariant ideals of $B(G)$ correspond to closed ideals of $\widehat{G}$, it follows that every weak$^*$-closed $G$-invariant ideal of $B(G)$ is of the form $B_{S_p}(G)$ for some $p \in [2,\infty]$.

    Conversely, for each $p \in [2,\infty]$, since $S_p$ is a closed ideal of $\widehat{G}$, the corresponding Fourier--Stieltjes space $B_{S_p}(G)$ is a weak$^*$-closed $G$-invariant ideal of $B(G)$. Hence, the weak$^*$-closed $G$-invariant ideals of $B(G)$ are precisely the Fourier--Stieltjes spaces $B_{S_p}(G)$, for all $p \in [2,\infty]$.

\end{proof}

\begin{thm}
    Let $G = \mathcal{K}^2 \rtimes H$. Then, for each $p\in[2,\infty]$, we have
    \begin{align*}
        B_{S_p}(G) = B_r(G) + B_{L^p}(H) \circ q.
    \end{align*}
\end{thm}

\begin{proof}
    We choose a representation $\sigma_p$ of $H$ such that $B_{\sigma_p}(H) = B_{L^p}(H).$
    Then $\sigma_p \circ q$ is a representation of $G$, and $B_{\sigma_p \circ q}(G) = B_{\sigma_p}(H) \circ q = B_{L^p}(H) \circ q$ (see \cite[Corollary 2.2.4]{Kan-Lau}). Also, $B_{\lambda_G}(G) = B_r(G).$ Since $S_p = \widehat{G}_r \cup (\widehat{H}_{L^p} \circ q)$, the representation $\rho := \lambda_G \oplus (\sigma_p \circ q)$ corresponds to $S_p$. Hence, $B_{S_p}(G) = B_\rho(G)$, and it follows from \cite[THEOREME 3.9 (iii)]{Arsac} that $B_\rho(G) = B_{\lambda_G \oplus (\sigma_p \circ q)}(G) = B_{\lambda_G}(G)+B_{\sigma_p \circ q}(G).$ Thus,
    \begin{align*}
        B_{S_p}(G) = B_r(G) + B_{L^p}(H) \circ q.
    \end{align*}
\end{proof}

    We have already shown that, for $H = \text{SL}(2,\mathcal{K})$, where
    $\mathcal{K} = \mathbb{R}$, $\mathbb{C}$, or $\mathbb{Q}_p$, the two classes of $L^p$-Fourier--Stieltjes algebras coincide. Thus, throughout the preceding arguments, we have removed the labels ``BG`` and ``KLQ`` and simply written $B_{L^p}(H)$.

\begin{prop}
    Let $G = \mathcal{K}^{3}\rtimes H$, where $H = \mathcal{K}^{2} \rtimes \text{SL}(2,\mathcal{K})$, and $q_H: H \to \text{SL}(2,\mathcal{K})$
    be the canonical quotient map. Then the closed ideals of $\widehat{G}$ are precisely the sets
    \begin{align*}
        J_p = \widehat{G}_r \cup \bigcup_{k \in \mathcal{K}}
        \left\{
        \chi_{(0,k)} \times (\sigma\circ q_H):
        \sigma \in \widehat{\text{SL}(2,\mathcal{K})}_{L^p}
        \setminus \widehat{\text{SL}(2,\mathcal{K})}_r
        \right\},
        \hspace{3mm} p \in [2,\infty].
    \end{align*}

\end{prop}

\begin{proof}
By Example~\ref{eg:3.14},
\begin{align*}
\widehat{G} \setminus \widehat{G}_r = \bigcup_{k\in\mathcal{K}}
\left\{
\chi_{(0,k)} \times \pi:
\pi \in \widehat{H}\setminus \widehat{H}_r
\right\}.
\end{align*}
Hence, by Theorem~\ref{thm:3.19}, every closed ideal of $\widehat{G}$ is of the form
\begin{align*}
    \widehat{G}_r \cup \bigcup_{k\in\mathcal{K}} \left\{
\chi_{(0,k)} \times \pi: \pi \in S \setminus \widehat{H}_r \right\},
\end{align*}
where $S$ is a closed ideal of $\widehat{H}$.
Since $H = \mathcal{K}^{2} \rtimes \text{SL}(2,\mathcal{K})$, Proposition~\ref{cor:3.20} shows that the closed ideals of $\widehat{H}$ are precisely
\begin{align*}
S_p = \widehat{H}_r \cup \left\{ \sigma \circ q_H: \sigma \in \widehat{\text{SL} (2,\mathcal{K})}_{L^p} \setminus \widehat{\text{SL} (2,\mathcal{K})}_r \right\},
\hspace{3mm} p \in [2,\infty].
\end{align*}
Then, 
\begin{align*}
S_p \setminus \widehat{H}_r = \left\{ \sigma \circ q_H:
\sigma \in \widehat{\text{SL}(2,\mathcal{K})}_{L^p}
\setminus \widehat{\text{SL}(2,\mathcal{K})}_r \right\}.
\end{align*}
This completes the proof.
\end{proof}

\section*{Acknowledgment}
The author would like to thank his Ph.D. supervisor, Matthew Wiersma, for introducing him to this project and for helpful discussions and suggestions throughout this research.

%%%%%%%%%%%%%%%%%%%%%%%%%%%%%%%%%%%%%%%%%%%%%%%%%%%%%%%%%%%%%%%%%%%%%%%%%%%%%%%%%%%%%%%%%%%%%%
\bibliographystyle{unsrtnat}

\begin{thebibliography}{23}
\providecommand{\natexlab}[1]{#1}
\providecommand{\url}[1]{\texttt{#1}}
\expandafter\ifx\csname urlstyle\endcsname\relax
  \providecommand{\doi}[1]{doi: #1}\else
  \providecommand{\doi}{doi: \begingroup \urlstyle{rm}\Url}\fi

\bibitem{Arsac}
G. Arsac,
\newblock{\it Sur l'espace de Banach engendr\'e par les coefficients d'une repr\'esentation unitaire},
\newblock Publ. D\'ep. Math. (Lyon) 13 (1976), no. 2, 1-101.

\bibitem{Bekka}
B. Bekka and P. de la Harpe,
\newblock{\it Unitary representations of groups, duals, and characters},
\newblock Mathematical Surveys and Monographs, vol. 250, Amer. Math. Soc., (2020).

\bibitem{Bek-Ech}
B. Bekka and S. Echterhoff,
\newblock{\it On unitary representations of algebraic groups over local fields},
\newblock J. Represent. Theory., Vol. 25, (2021) 508-526.


\bibitem{Br-Gn} 
N.P. Brown and E.P. Guentner, 
\newblock {\it   New $C^*$-completions of discrete groups and related spaces}, 
\newblock  Bull. Lond. Math. Soc. 45 (2013), no. 6, 1181-1193.

\bibitem{Carter}
M. Carter,
\newblock{\it On the unitary representation theory of locally compact contraction groups},
\newblock J. Lie Theory 34 (2024), no. 4, 911--956.

\bibitem{de Laat}
T. de Laat and T. Siebenand.
\newblock{\it Exotic group C$^*$-algebras of simple lie groups with real rank one},
\newblock Ann. Inst. Fourier (Grenoble) 71 (2021), no. 5, 2117-36.


\bibitem{FTW}
B. Forrest, Z. Tanko and M. Wiersma,
\newblock{\it Exotic ideals in the Fourier--Stieltjes algebra of a locally compact group},
\newblock J. Math. Anal. Appl. 506 (2022), Paper No. 125528.

\bibitem{HLS-Tr-IMRN}
D. Heinig, T. de Laat and T. Siebenand,
\newblock{\it Group $C^*$-algebras of locally compact groups acting on trees},
\newblock Int. Math. Res. Not. IMRN (2024), no. 10, 8520-8539. 



\bibitem{KLQ-Ex-NY}
S. Kaliszewski, M.B. Landstad and J. Quigg,
\newblock{\it Exotic group $C^*$-algebras in noncommutative duality},
\newblock New York J. Math. 19 (2013) 689-711.

\bibitem{KLQ-Oi-NY}
S. Kaliszewski, M.B. Landstad and J. Quigg,
\newblock{\it Ordered invariant ideals of Fourier--Stieltjes algebras},
\newblock New York J. Math. 24 (2018) 1039-1055.

\bibitem{Kan-Lau}
E. Kaniuth, A. T. -M. Lau,
\newblock{\it Fourier and Fourier--Stieltjes algebras on locally compact groups},
\newblock Mathematical Surveys and Monographs, Vol. 231, Amer. Math. Soc. (2018).

\bibitem{Kan-Tay}
E. Kaniuth, K.F. Taylor,
\newblock{\it Induced representations of locally compact groups},
\newblock Cambridge Tracts in Mathematics, Vol. 197, Cambridge Univ. Press (2013).

\bibitem{Kechris}
A. S. Kechris,
\newblock{\it Classical descriptive set theory},
\newblock Graduate Text in Mathematics, Vol. 156, Springer-Verlag, New York (1995).


\bibitem{Koo-RLC}
T. H. Koornwinder (ed.),
\newblock{\it Representations of locally compact groups
with applications, Part II},
\newblock MC Syllabus 38.2, Mathematisch Centrum,
Amsterdam (1979).


\bibitem{SW-Ex-Fnc}
E. Samei and M. Wiersma,
\newblock{\it Exotic $C^*$-algebras of geometric groups},
\newblock J. Funct. Anal. 286 (2) (2024), Paper No. 110228.

\bibitem{Sie-FS-arx}
T. Siebenand,
\newblock{\it On the ideal structure of the Fourier--Stieltjes algebra of certain groups},
\newblock Preprint, arXiv: 2005.01772 (2020).



\bibitem{WRS-Fnc}
M. Wiersma,
\newblock{\it $L^p$-Fourier and Fourier--Stieltjes algebras for locally compact groups},
\newblock J. Funct. Anal. 269 (2015), no. 12, 3928-3951.

\end{thebibliography}

%%%%%%%%%%%%%%%%%%%%%%%%%%%%%%%%%%%%%%%%%%%%%%%%%%%%%%%%%%%%

\bigskip
\noindent
\textsc{Hamidreza Shirnejad}\\
Department of Mathematics, University of Manitoba, Winnipeg, Manitoba, Canada\\
\textit{E-mail address:} \texttt{shirnejh@myumanitoba.ca}

\end{document}